\documentclass[10pt,notitlepage]{article}
\usepackage{amsfonts}
\usepackage{amsmath}
\usepackage{epsfig}
\usepackage{graphicx}
\title{Simplicity of QC($\mathbb{S}^{n}$) and LIP($\mathbb{S}^{n}$)}
\author{James Giblin}
\newtheorem{thm}{Theorem}[section]
\newtheorem{cor}[thm]{Corollary}
\newtheorem{lem}[thm]{Lemma}

\newtheorem{defn}[thm]{Definition}
\newtheorem{rmk}[thm]{Remark}
\newtheorem{exa}[thm]{Example}

\begin{document}

\maketitle

\begin{center}
\begin{small}
\textsc{ABSTRACT\\}
\end{small}
\vspace{.3cm}
We will prove that the groups QC($\mathbb{S}^{n}$) and LIP($\mathbb{S}^{n}$) ($n\geq2$) are simple.

\end{center}

\vspace{.5cm}

\section{Introduction}

A group $G$ is simple if its only normal subgroups are itself and the trivial group. Let QC($\mathbb{S}^{n}$) and LIP($\mathbb{S}^{n}$) denote the orientation preserving quasiconformal and bilipschitz homeomorphisms of $\mathbb{S}^{n}$ respectively (We will implicitly assume throughout that the dimension $n$ is greater than 1). In this paper we prove the following statement,

\begin{thm}
The groups \textnormal{QC($\mathbb{S}^{n}$)} and \textnormal{LIP($\mathbb{S}^{n}$)} are simple.
\end{thm}

Our proof falls into two stages. First, we show that the subgroups
of QC($\mathbb{S}^{n}$) and LIP($\mathbb{S}^{n}$) which are
generated by the elements supported on open balls are simple. In
the second stage, we show that \emph{every} element of
QC($\mathbb{S}^{n}$) or LIP($\mathbb{S}^{n}$) can be written as a
composition of elements which are all supported on open balls. The
main part of the argument is Lemma \ref{lem:stab} which proves
that there are neighbourhoods of the identity in
QC($\mathbb{S}^{n}$) and LIP($\mathbb{S}^{n}$) which satisfy this
condition.

\vspace{.3cm}

Anderson \cite{AN} proved that if $G$ is a group of
homeomorphisms of a Hausdorff space $X$ and $K$ is a basis for
$X$ satisfying certain conditions, then the group $G_{0}$ generated by the homeomorphisms which are supported on elements of $K$ is simple. In fact, he showed that if
$h$ is any non trivial element of $G$, then every element of $G_{0}$
can be written as a product of four conjugates of $h$ and
$h^{-1}$. We demonstrate in detail that the group of
quasiconformal homeomorphisms of $\mathbb{S}^{n}$ or the group of Lipschitz
homeomorphisms of $\mathbb{S}^{n}$ and the basis consisting of open balls
satisfy the necessary conditions to apply his argument and
give the proof that these groups will therefore be simple.

\vspace{.3cm}

In his paper, Anderson also gave some conditions that the whole
group $G$ will be generated by $G_{0}$. Namely that $X$ can be
expressed as the union of two basis elements, and that $G$ has the
"Annulus property". Here, $G$ satisfies the annulus property if
for every $g\in G$ and $k\in K$ such that there exists $k'\in K$
with $k\cup g(k)\subset k'$, there exists $g_{0}\in G_{0}$ such
that $g$ and $g_{0}$ agree on $k$. Now, $\mathbb{S}^{n}$ can
obviously be written as the union of two open balls and in the
process of proving Lemma \ref{lem:stab} we show that there are
neighbourhoods of $id$ in the groups QC($\mathbb{S}^{n}$) and
LIP($\mathbb{S}^{n}$) which satisfy the annulus property. In these
cases, this turns out to be enough.

\vspace{.3cm}

A homeomorphism $h:X\to X$ is called \emph{stable} if it can be
written as a finite composition of homeomorphisms, each of which
is the identity on some open set. We say a group of homeomorphisms
$G$ is \emph{stable} if every element of $G$ is (See \cite{BG} and \cite{BR} for more information on stable homeomorphisms). Kirby \cite{KI}
showed that the group of homeomorphisms of the torus $\mathbb{T}^{n}$ is
stable. He used this, along with an immersion of $\mathbb{T}^{n}$ minus a
point into $\mathbb{R}^{n}$ to show that the group of
homeomorphisms of $\mathbb{R}^{n}$ is locally contractible.
Edwards and Kirby \cite{EK} then went on to extend these results
to arbitrary manifolds, in particular, showing that the
homeomorphism group of a compact manifold $M$ is locally
contractible.

\vspace{.3cm}

The stable homeomorphism conjecture says that the orientation
preserving homeomorphisms of $\mathbb{R}^{n}$ are stable. This is
known to be true for $n\neq 4$ \cite{KI} whereas for $n=4$ this
question is unsolved. A simple application of Anderson's results
concludes from the stable homeomorphism conjecture that the group
of orientation preserving homeomorphisms of $\mathbb{S}^{n}$ is simple for
$n\neq 4$.

\vspace{.3cm}

In general, if $M$ is a non compact manifold then its
group of homeomorphisms will not be simple, as the subgroup
consisting of elements with compact support forms a non-trivial
normal subgroup. Likewise for manifolds with boundary, the
subgroup consisting of elements which act as the identity in some
neighbourhood of the boundary will also be a non-trivial normal
subgroup. In the case of an arbitrary closed manifold the homeomorphism group
has generally got more than one path component, with the component
of the identity (which consists of all the stable homeomorphisms)
forming a non-trivial normal subgroup. It is still an interesting question to ask whether these subgroups are themselves simple.

\vspace{.3cm}

For the case of diffeomorphisms, let $M$ be an $n$-dimensional smooth manifold. Define Diff$(M,r)$ to be the group of $C^{r}$ diffeomorphisms of $M$ which are isotopic to the identity through compactly supported $C^{r}$ isotopies (An isotopy has compact support if there is compact set outside of which the isotopy acts trivially). In \cite{M1} and \cite{M2} Mather proved that if $\infty\geq r\neq n+1$ then Diff$(M,r)$ is perfect (equal to its own commutator subgroup) with the case $r=\infty$ being proved by Thurston (See \cite{T1}). Now, in \cite{EP}, Epstein showed that for certain groups of homeomorphisms, the commutator subgroup is simple. In particular, Diff$(M,r)$ is simple when $\infty\geq r\neq n+1$.

\vspace{.3cm}

Simplicity questions have also been addressed in the symplectic and volume preserving categories. Perhaps most notably Banyaga \cite{BA} proved that for a closed symplectic manifold $M$ the group of Hamiltonian symplectomorphisms Ham$(M,\omega)$ is perfect and then used Epstein's argument, which also applies in the symplectic case, to deduce that it is simple. On a non compact manifold, if one considers the component of the identity in the group of compactly supported symplectomorphisms $\mathrm{Symp}_{0}^{c}(M,\omega)$, then there exists a surjective homomorphism from $\mathrm{Symp}_{0}^{c}(M,\omega)$ to $\mathbb{R}$. Banyaga proved that the kernel of this homomorphism is simple and also showed an analogous result in the volume preserving case. In general however, little is known about the normal subgroups of $\mathrm{Symp}(M,\omega)$. For example we currently do not know whether the group of symplectomorphisms of $\mathbb{R}^{2n}$ with its standard symplectic form is perfect. A stark contrast to the volume preserving case, where Mcduff \cite{MC} has shown that $\mathrm{Diff_{vol}}(\mathbb{R}^{n})$ is perfect for every dimension $n\geq 3$.

\vspace{.3cm}

In \cite{SU} Sullivan proved the stable homeomorphism
conjecture and the annulus conjecture in the locally quasiconformal
and Lipschitz categories in all dimensions. We go through his proof of the stable homeomorphism conjecture in detail, and use the result to show simplicity of QC($\mathbb{S}^{n}$) and
LIP($\mathbb{S}^{n}$).

\vspace{.3cm}

I would like to thank Vladimir Markovic for introducing me to this problem and for many helpful discussions regarding this work.

\section{QC and LIP homeomorphisms supported on open balls generate a simple group.}

In this section we start by introducing some notation and terminology and go on to prove that every element of QC($\mathbb{S}^{n}$) or LIP($\mathbb{S}^{n}$) which is supported on an open ball can be written as a composition of four conjugates of an arbitrary non-trivial element and its inverse. More importantly for us, this gives us the easy corollary that the subgroups of QC($\mathbb{S}^{n}$) and LIP($\mathbb{S}^{n}$) generated by those elements supported on open balls are simple.

\vspace{.3cm}

Let $d_{\mathbb{R}^{n}}$ denote the usual Euclidean distance on
$\mathbb{R}^{n}$, we set
$$
\mathbb{S}^{n-1}=\{x\in\mathbb{R}^{n}:d_{\mathbb{R}^{n}}(x,0)=1\}
$$
to be the unit sphere in $\mathbb{R}^{n}$ and denote the standard spherical metric on $\mathbb{S}^{n}$ by $d_{\mathbb{S}^{n}}$. Let
$$
\mathbb{B}^{n}=\{x\in\mathbb{R}^{n}:d_{\mathbb{R}^{n}}(x,0)<1\}
$$
denote the unit ball in $\mathbb{R}^{n}$, and let $d_{\mathbb{B}^{n}}$ denote its hyperbolic metric. We now define,
$$
B_{X}(c,\epsilon)=\{x\in X:d_{X}(x,c)<\epsilon\}
$$
to be the open ball of radius $\epsilon>0$ around $c\in X$, where $X$ is one of $\mathbb{R}^{n}$, $\mathbb{S}^{n}$ or $\mathbb{B}^{n}$ and $d_{X}$ is the appropriate metric.

We also define a metric on the space of maps from $\mathbb{S}^{n}$ to itself by,
$$
\tilde{d}(f,g)=sup\{d(f(x),g(x)):x\in\mathbb{S}^{n}\}.
$$

The group of M\"obius transformations M\"ob($\mathbb{S}^{n}$) is the transformation group of $\mathbb{S}^{n}$ generated by inversions in ($n-1$)-spheres and M\"ob$^{+}(\mathbb{S}^{n})$ is the index two subgroup consisting of those which preserve orientation. By identifying $\mathbb{S}^{n}$ with $\mathbb{R}^{n}\cup\{\infty\}$ we can think of $\mathbb{B}^{n}$ as a subset of $\mathbb{S}^{n}$. If we define M\"ob($\mathbb{B}^{n}$)=$\{f\in$M\"ob$(\mathbb{S}^{n}):f(\mathbb{B}^{n})=\mathbb{B}^{n}\}$ then M\"ob($\mathbb{B}^{n}$) is generated by inversions in ($n-1$)-spheres orthogonal to $\mathbb{S}^{n-1}=\partial \mathbb{B}^{n}$ and is the isometry group of $\mathbb{B}^{n}$ with its hyperbolic metric. As before
M\"ob$^{+}(\mathbb{B}^{n})\subset$M\"ob$(\mathbb{B}^{n})$ denotes those M\"obius
transformations which preserve orientation. The groups
M\"ob$(\mathbb{B}^{n})$ and M\"ob$^{+}(\mathbb{B}^{n})$ preserve $\mathbb{S}^{n-1}$ and in the obvious way are isomorphic to M\"ob$(\mathbb{S}^{n-1})$ and
M\"ob$^{+}(\mathbb{S}^{n-1})$ respectively.

\vspace{.3cm}

\begin{defn}
Let $(X,d_{X})$ and $(Y,d_{Y})$ be metric spaces and let $f:X\to Y$ be a map between them. Then if there is a constant $L>0$ such that
$$
d_{Y}(f(x),f(y))\leq L d_{X}(x,y)
$$
for every $x,y\in X$, then $f$ is $L$-\emph{Lipschitz} . Furthermore, if we also have that
$$
d_{Y}(f(x),f(y))\geq d_{X}(x,y)/L
$$
for every $x,y\in X$ then $f$ is $L$-\emph{bilipschitz}.
\end{defn}

\begin{defn}
If $f:(X,d_{X})\to (Y,d_{Y})$ is such that every $x\in X$ has a neighbourhood $U$ such that the restriction of $f$
to $U$ is Lipschitz, then $f$ is \emph{locally Lipschitz} (similarly for the definition of locally
bilipschitz)
\end{defn}

\begin{rmk}
\textnormal{Metric spaces and locally Lipschitz maps form a category usually denoted by LIP. For this reason we will abbreviate locally Lipschitz to LIP. Moreover, if $f:(X,d_{X})\to (Y,d_{Y})$ is bijective and both
$f$ and $f^{-1}$ are LIP maps, then we say $f$ is a
\emph{lipeomorphism}.}
\end{rmk}

We write LIP($X$) to denote the space of
lipeomorphisms from $X$ to itself. Since every LIP map with compact
domain is Lipschitz, the space LIP($\mathbb{S}^{n}$) consists
of all the bilipschitz homeomorphisms from $\mathbb{S}^{n}$ to itself (the point here is that we don't have to say \emph{locally} bilipschitz). It is a standard result
that every diffeomorphism from $\mathbb{S}^{n}$ to itself lies in
LIP($\mathbb{S}^{n}$).

\begin{rmk}
 We will only consider those homeomorphisms which preserve the orientation of $\mathbb{S}^{n}$.
\end{rmk}

Let $\Omega,\Omega'$ be subdomains of $\mathbb{R}^{n}$ and $f:\Omega\to\Omega'$ be a homeomorphism. Fix $x_{0}\in\Omega$ and for small enough $\epsilon>0$ consider the quantities
$$
D_{l}(\epsilon)=inf\{d_{\mathbb{R}^{n}}(f(x_{0}),y):y\in f(\partial B_{\mathbb{R}^{n}}(x_{0},\epsilon))\}
$$
and
$$
D_{u}(\epsilon)=sup\{d_{\mathbb{R}^{n}}(f(x_{0}),y):y\in f(\partial B_{\mathbb{R}^{n}}(x_{0},\epsilon))\}
$$
Let $K_{x_{0}}(\epsilon)=D_{u}(\epsilon)/D_{l}(\epsilon)\geq 1$. If there exists $K<\infty$ such that $K_{x_{0}}=limsup_{\epsilon\to 0}K_{x_{0}}(\epsilon)<K$ for every $x_{0}\in\Omega$ then we say $f$ is $K$-quasiconformal. We say $f$ is quasiconformal if $f$ is $K$-quasiconformal for some $K$. If every point $x_{0}\in\Omega$ has a neighbourhood $U$ such that $f|_{U}$ is $K$-quasiconformal, then we say $f$ is locally $K$-quasiconformal.

\vspace{.3cm}

Locally quasiconformal homeomorphisms between open subsets of $\mathbb{R}^{n}$ form a pseudogroup on $\mathbb{R}^{n}$ (See \cite{TH} for the definition of a pseudogroup and more detail). In the usual way, this allows us to construct \emph{quasiconformal manifolds}. These are manifolds whose transition maps are locally quasiconformal. Furthermore, if we have a homeomorphism $f:M\to N$ between two such manifolds, then we say it is locally quasiconformal if the following holds. For any pair of charts $\phi:M\to\mathbb{R}^{n}$ and $\psi:N\to\mathbb{R}^{n}$ the composition $\psi\circ f\circ\phi^{-1}$ is locally quasiconformal where defined. In particular, we let QC($\mathbb{S}^{n}$) denote the group of orientation preserving quasiconformal mappings from $\mathbb{S}^{n}$ to itself. Notice here that the standard smooth structure on $\mathbb{S}^{n}$ automatically induces a locally quasiconformal structure and as in the case of LIP($\mathbb{S}^{n}$), QC($\mathbb{S}^{n}$) contains all the orientation preserving diffeomorphisms of $\mathbb{S}^{n}$.

\begin{rmk}
\textnormal{The space of lipeomorphisms between open subsets of $\mathbb{R}^{n}$ also forms a pseudogroup. So we could equally well have defined LIP($\mathbb{S}^{n}$) to be the group of orientation preserving bilipschitz homeomorphisms of $\mathbb{S}^{n}$ where its Lipschitz structure is induced from the standard smooth structure.}
\end{rmk}

Throughout the rest of the paper we take G($\mathbb{S}^{n}$) to denote either \textnormal{LIP($\mathbb{S}^{n}$)} or \textnormal{QC($\mathbb{S}^{n}$)}.

\vspace{.3cm}

Take $\mathcal{B}$ to be the basis for
the induced topology on $\mathbb{S}^{n}$ consisting of the open balls
$\{B_{\mathbb{S}^{n}}(c,\epsilon): c\in \mathbb{S}^{n},~\epsilon\in(0,\pi) \}$. This
ensures that for all $B\in \mathcal{B}$, $\mathbb{S}^{n}\setminus B$ is
homeomorphic to $\bar{\mathbb{B}}^{n}$ by excluding the possibility that
$B=\mathbb{S}^{n}$ or $\mathbb{S}^{n}$ minus a point. We do this to ensure that
given two elements $B_{1},B_{2}\in\mathcal{B}$ we can find $g\in
G(\mathbb{S}^{n})$ such that $g(B_{1})=B_{2}$.

\vspace{.3cm}

If $g\in\mathrm{G}(\mathbb{S}^{n})$
is not the identity, then we say $g$ is \emph{supported} on $B\in
\mathcal{B}$ if $g$ is the identity outside $B$. Consequently, if
$B_{1}\subset B_{2}$ and $g\in \mathrm{G}(\mathbb{S}^{n})$ is supported on
$B_{1}$ then $g$ is also supported on $B_{2}$. We set
$\mathrm{G}_{0}(\mathbb{S}^{n})$ to be those elements of
$\mathrm{G}(\mathbb{S}^{n})$ supported on elements of $\mathcal{B}$.

\vspace{.3cm}

We have already mentioned that for all $B_{1},B_{2}\in
\mathcal{B}$ there exists $g\in \mathrm{G}(\mathbb{S}^{n})$ with
$g(B_{1})=B_{2}$. In particular, if $B_{1},B_{2}\in \mathcal{B}$
and $g\in \mathrm{G}_{0}(\mathbb{S}^{n})$ is an element supported on
$B_{1}$, then $g$ is conjugate to an element of
$\mathrm{G}_{0}(\mathbb{S}^{n})$ supported on $B_{2}$. Furthermore, we can
choose $g$ to lie in M\"ob$^{+}(\mathbb{S}^{n})\subset\mathrm{G}(\mathbb{S}^{n})$ so that if
$B\in\mathcal{B}$ then
$g(B)$ will also lie in $\mathcal{B}$.

\begin{defn}
A pair $(\{B_{i}\}_{i \geq 0},\rho)$ consisting of a sequence
$\{B_{i}\}_{i \geq 0}$ of elements of $\mathcal{B}$, and $\rho\in
\mathrm{G}_{0}(\mathbb{S}^{n})$ will be called a
\textbf{$\sigma$-sequence} if,
\begin{enumerate}
\item{$B_{i}\cap B_{j}=\phi$ if $i\neq j$}
\item{There exists $x_{0}\in \mathbb{S}^{n}$ such that every neighbourhood
of $x_{0}$ contains all but finitely many of the $B_{i}$}
\item{$\rho(B_{i})=B_{i+1}$}
\end{enumerate}
\end{defn}

This means that if $\rho$ is supported on $B\in \mathcal{B}$ then
$\bigcup_{i\geq 0}B_{i}\subset B$ (since condition 1 means that
$\rho$ cannot be the identity anywhere in $\bigcup_{i\geq
0}B_{i}$). If $(\{B_{i}\}_{i \geq 0},\rho)$ is a sigma sequence,
we will say it is supported on $B\in\mathcal{B}$ if $\rho$ is.

\vspace{.3cm}

\begin{exa}\label{exa:seq}
\textnormal{We now use the fact that M\"ob$^{+}(\mathbb{S}^{n})\subset
\mathrm{G}(\mathbb{S}^{n})$ to construct a $\sigma$-sequence explicitly.
Consider $\mathbb{S}^{n}$ as $\mathbb{R}^{n}\cup\{\infty\}$ then the
M\"obius transformation,}
\end{exa}

$$
f(\textbf{x})=\frac{\textbf{x}}{r^{2}}, f(\infty)=\infty
$$
with $r\in\mathbb{R}\setminus\{0\}$ is the time $r^{2}$ flow of
the complete vector field,
$$
\mathcal{X}(x_{1},\dots,x_{n})=-\Sigma_{i} x_{i}\frac{\partial}{\partial
x_{i}}
$$

Let $\varphi:\mathbb{S}^{n}\to\mathbb{R}$ be a smooth function
which is identically one on $B_{\mathbb{S}^{n}}(0,1)$ and zero
outside $B_{\mathbb{S}^{n}}(0,1+\epsilon)$ for some $0<\epsilon<1$
and take $\rho$ to be the time $r^{2}$ flow of
$\varphi\cdot\mathcal{X}$. Then $\rho$ will be a smooth map which
is the same as $f$ on $B_{\mathbb{S}^{n}}(0,1)$ and the identity
outside $B_{\mathbb{S}^{n}}(0,1+\epsilon)$ and so is in
$\mathrm{G}_{0}(\mathbb{S}^{n})$. Now take $B\in \mathcal{B}$ with
$B\subset B_{\mathbb{S}^{n}}(0,1)$ and $\rho(B)\cap B=\phi$ and
define $\{B_{i}\}=\rho^{i}(B)$ $i\geq 0$, these will all be in
$\mathcal{B}$ since $\rho$ is conformal on
$B_{\mathbb{S}^{n}}(0,1)$ and $(\{B_{i}\},\rho)$ will be a
$\sigma$-sequence supported on $B_{\mathbb{S}^{n}}(0,1+\epsilon)$
with $x_{0}=0$. Take any $B\in \mathcal{B}$, then since there
exists $g\in \mathrm{G}(\mathbb{S}^{n})$ with
$g(B_{\mathbb{S}^{n}}(0,1))=B$, we have that
$(\{g(B_{i})\},g\circ\rho\circ g^{-1})$ is a $\sigma$-sequence
supported on $g(B_{\mathbb{S}^{n}}(0,1+\epsilon))$.

\vspace{.3cm}

Let $(\{B_{i}\},\rho)$ be a $\sigma$-sequence supported on $B\in
\mathcal{B}$, for each $i$ take $g_{i}\in \mathrm{G}_{0}(\mathbb{S}^{n})$
which is supported on $B_{i}$. We would like to be able to say that $\prod_{i\geq 0}g_{i}=\cdots\circ g_{1}\circ g_{0}$ is also in $\mathrm{G}_{0}(\mathbb{S}^{n})$. If $g_{i}\in \mathrm{QC}_{0}(\mathbb{S}^{n})$ for each $i$
and there exists K $\geq 1$ such that
the $g_{i}$ are all K-quasiconformal, then the composition
$\prod_{i\geq 0}g_{i}$ is also K-quasiconformal and hence in
$\mathrm{QC}_{0}(\mathbb{S}^{n})$. In particular, if $\rho|_{\cup B_{i}}$
is conformal (this was the case in the example we constructed
above) and we take $g_{0}$ supported on $B_{0}$ which is
K-quasiconformal. Then the conjugates $\rho^{-i}\circ g_{0}\circ \rho^{i}$
will also be K-quasiconformal and supported on $B_{i}$. Hence,
 the composition $\prod_{i\geq 0}\rho^{-i}\circ g_{0}\circ \rho^{i}$ will also be K-quasiconformal,
and thus lie in $\mathrm{QC}_{0}(\mathbb{S}^{n})$.

\vspace{.3cm}

Similarly, if for each $i$, $g_{i}\in \mathrm{LIP}_{0}(\mathbb{S}^{n})$ is supported on $B_{i}$
and there exists $L\geq 1$ such that
the Lipschitz constants of the $g_{i}$ are bounded by $L$, then the composition
$\prod_{i\geq 0}g_{i}$ is $L$-Lipschitz and hence in
$\mathrm{LIP}_{0}(\mathbb{S}^{n})$. In particular, if we take $g_{0}$ supported on $B_{0}$ which is
$L$-Lipschitz. Then again using $\rho$ from the example above, the conjugates $\rho^{-i}\circ g_{0}\circ \rho^{i}$
will also be $L$-Lipschitz and supported on $B_{i}$. Hence,
 the composition $\prod_{i\geq 0}\rho^{-i}\circ g_{0}\circ \rho^{i}$ will also be $L$-Lipschitz,
and thus lie in $\mathrm{LIP}_{0}(\mathbb{S}^{n})$. We will use these observations to prove the following theorem.

\begin{thm}\label{thm:conj}
Given a non-trivial $h\in \mathrm{G}(\mathbb{S}^{n})$ , any element $g\in
\mathrm{G}_{0}(\mathbb{S}^{n})$ is the product of four conjugates of $h$
and $h^{-1}$.
\end{thm}

\noindent\textsc{Proof.} As $h$ is not the identity, we can find
$B\in \mathcal{B}$ with $h^{-1}(B)\cap B=\phi$. Let
$(\{B_{i}\},\rho)$ denote a $\sigma$-sequence with $\cup_{i} B_{i}\subset B$ obtained as in Example \ref{exa:seq}, so that $\rho|_{B}$ is conformal (Note that $\rho$ will be supported on some larger ball containing $B$). Since every element of
$\mathrm{G}_{0}(\mathbb{S}^{n})$ can be conjugated in $\mathrm{G}(\mathbb{S}^{n})$
to one supported on $B_{0}$ it suffices to show that any $g$
supported on $B_{0}$ is the product of four conjugates of $h$ and
$h^{-1}$.

\vspace{.3cm}

So let $g\in \mathrm{G}_{0}(\mathbb{S}^{n})$ be supported on $B_{0}$, then
the map,
$$
f=\prod_{n\geq 0}\rho^{n}\circ g\circ \rho^{-n}
$$
is an element of $\mathrm{G}_{0}(\mathbb{S}^{n})$ supported on $\cup_{i\geq 0}B_{i}\subset
B$. Since the conjugate \mbox{$h^{-1}\circ  f^{-1}\circ  h$} is then supported on
$h^{-1}(B)$ and $h^{-1}(B)\cap B=\phi$, we can think of the
commutator,
$$
f_{1}:=h^{-1}\circ f^{-1}\circ  h\circ  f =(h^{-1}\circ  f^{-1}\circ  h)\circ  f
$$
as two actions on disjoint sets.

\vspace{.3cm}

The same can be said about,
$$
f_{2}:=\rho\circ  h^{-1}\circ  f\circ  h\circ  f^{-1}\circ  \rho^{-1}=(\rho\circ  h^{-1}\circ  f\circ  h\circ
\rho^{-1})\circ  (\rho\circ  f^{-1}\circ  \rho^{-1})
$$
here $\rho\circ  f^{-1}\circ  \rho^{-1}$ is supported on $B$ whereas $\rho\circ
h^{-1}\circ  f\circ  h\circ  \rho^{-1}$ is supported on $h^{-1}(B)$.

\vspace{.3cm}

So $f_{2}\circ f_{1}$ can be examined in terms of its action on
$B$ and $h^{-1}(B)$ respectively. Doing this, one can easily see
that $f_{2}\circ f_{1}$ is the identity on $h^{-1}(B)$ and acts as
$g$ on $B$. Since $g$ is supported on $B_{0}\subset B$,
$g=f_{2}\circ f_{1}$, but we can write,
$$
g=f_{2}\circ f_{1}=(\rho\circ  h^{-1}\circ \rho^{-1})\circ (\rho\circ  f\circ  h\circ  f^{-1}\circ
\rho^{-1})\circ  (h^{-1})\circ  (f^{-1} \circ h \circ f)
$$
as a product of four conjugates of $h$ and $h^{-1}$.
\marginpar{$\square$}

\begin{cor}
The subgroup of $\mathrm{G}(\mathbb{S}^{n})$ generated by
$\mathrm{G}_{0}(\mathbb{S}^{n})$ is simple.
\end{cor}

\noindent\textsc{Proof.} Denote by $\mathrm{G}(\mathrm{G}_{0}(\mathbb{S}^{n}))$ the subgroup of
$\mathrm{G}(\mathbb{S}^{n})$ generated by $\mathrm{G}_{0}(\mathbb{S}^{n})$ and let $H$ be a non-trivial normal subgroup. We need to show that $H=\mathrm{G}(\mathrm{G}_{0}(\mathbb{S}^{n}))$. So let $g\in\mathrm{G}(\mathrm{G}_{0}(\mathbb{S}^{n}))$ and take $h\in H$ which is not the identity. We can write $g$ as a finite composition $g_{n}\circ\cdots\circ g_{1}$ where each $g_{i}\in\mathrm{G}_{0}$ for $i=1\dots n$. Then using Theorem \ref{thm:conj} we can write each $g_{i}$ as a product of four conjugates of $h$ and $h^{-1}$, but since $H$ is normal this implies that each $g_{i}\in H$ and thus $g\in H$ which gives us the result.

\marginpar{$\square$}

\vspace{.3cm}

It now remains to show that $\mathrm{G}(\mathbb{S}^{n})$ is equal to
$G(\mathrm{G}_{0}(\mathbb{S}^{n}))$. We show this in the following
section.

\section{Sullivan groups}

In this section we introduce the concept of a Sullivan group, this
is a discrete group of hyperbolic isometries which gives rise to a
hyperbolic n-manifold with the property that the complement of any
point can be LIP-immersed into $\mathbb{R}^{n}$. We use this along
with two versions of the famous Schoenflies Theorem in the LIP
category to prove the main result, Lemma \ref{lem:stab}, which
states that there is a neighbourhood of the identity in
G($\mathbb{S}^{n}$) which is a subset of
$G(\mathrm{G}_{0}(\mathbb{S}^{n}))$. Using this, we prove that
$G(\mathrm{G}_{0}(\mathbb{S}^{n}))$ is equal to
G($\mathbb{S}^{n}$) which then must be simple.

\vspace{.3cm}

If $(X,d_{X})$ and $(Y,d_{Y})$ are metric spaces and $f:(X,d_{X})\to (Y,d_{Y})$ is a lipeomorphism onto its image (and hence injective)
then we say $f$ is a \emph{LIP embedding} of $X$ in $Y$. If every $x\in X$ has a
neighbourhood $U$ such that \hbox{$f|_{U}:U\to Y$} is a LIP embedding then we say $f$
is a \emph{LIP immersion}.

\vspace{.3cm}

Let $\Gamma\subset$M\"ob$^{+}(\mathbb{B}^{n})$ be a discrete group of
M\"obius transformations which acts freely, then the quotient
space $Q=\mathbb{B}^{n}/\Gamma$ is a smooth Hausdorff manifold and the
natural map $\pi:\mathbb{B}^{n}\to Q$ is a covering map. Note that since
$\Gamma$ acts properly and $\mathbb{B}^{n}$ is locally compact, the
property of $\Gamma$ being discrete is equivalent to it acting
properly discontinuously. $Q$ inherits a hyperbolic metric from
$\mathbb{B}^{n}$ defined in the usual way so that $\pi$ becomes a local
isometry and we write $B_{Q}(q,\epsilon)$ to denote the open ball of radius
$\epsilon$ around $q\in Q$ in this metric.

\begin{defn}\label{def:sul} A group $\Gamma$ of M\"obius transformations acting on $\mathbb{B}^{n}$ is a \textbf{Sullivan Group} if it satisfies the following;
\begin{enumerate}
\item{$\Gamma$ is discrete}
\item{$\Gamma$ acts freely}
\item{$Q=\mathbb{B}^{n}/\Gamma$ is closed (compact and without boundary)}
\item{For every $q\in Q$ there exists a LIP immersion $\iota:Q\setminus q\to \mathbb{R}^{n}$}
\end{enumerate}
\end{defn}

If $\Gamma$ is a Sullivan group we call the quotient space
$Q=\mathbb{B}^{n}/\Gamma$ a \emph{Sullivan manifold}. We have the following theorem from \cite{SU}.

\begin{thm}\label{thm:sul}
There exists a Sullivan manifold of every dimension $n\geq 2$.
\end{thm}

 Sullivan's proof of this theorem is deep and we take it on
faith here. As an easy example though, in dimension 2 any closed
surface of genus $\mathbf{g}\geq 2$ can be constructed as a
quotient of $\mathbb{B}^{n}$ by a group of M\"obius transformations
$\Gamma$ which is discrete and acts freely. We claim that all of
these groups are Sullivan groups.

\vspace{.3cm}

\begin{figure}[h]
\begin{center}
\includegraphics[width=5cm]{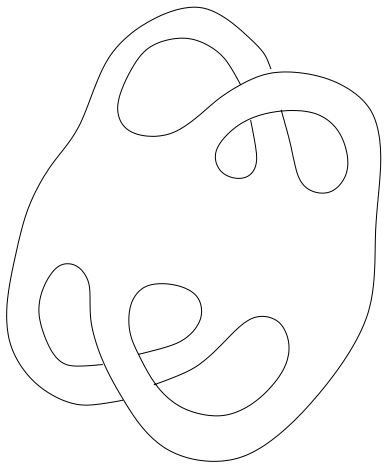}
\caption{Immersion of genus 2 surface minus a disc into
$\mathbb{R}^{2}$} \label{fig:immers}
\end{center}
\end{figure}

To see this it remains to show condition 4, that for every $q\in
Q$ there exists a LIP immersion $\iota:Q\setminus q\to
\mathbb{R}^{n}$. So take $q\in Q$ and remove the ball
$B_{Q}(q,\epsilon)$. The resulting surface is lipeomorphic
(diffeomorphic in fact) to a disc with $2\mathbf{g}$ strips
attached (see Figure \ref{fig:immers}), which can then be
LIP-immersed into $\mathbb{R}^{n}$. Letting $\epsilon$ tend to 0
then gives the required map (note however that the bilipschitz
constants will become arbitrarily large near $q$).

\vspace{.3cm}

\vspace{.3cm}

We will also need the following
two theorems which we state without proof.

\begin{thm}\label{thm:lip}
\textnormal{\textbf{(Lipschitz Schoenflies Theorem \cite{LV})}} Suppose that
$0<a<1$ and that $f$ is a LIP embedding of the annulus
$A=\bar{B}^{n}\setminus B_{\mathbb{R}^{n}}(0,a)$ into $\mathbb{S}^{n}$. Then
$f|_{S^{n-1}}$ can be extended to a LIP embedding
$f^{*}:\bar{B}^{n}\to \mathbb{S}^{n}$.
\end{thm}

\begin{thm}\label{thm:quan}
\textnormal{\textbf{(Quantitative Schoenflies Theorem \cite{TV})}} Let
$A=\bar{B}^{n}\setminus B_{\mathbb{R}^{n}}(0,1/2)$,
$\mathcal{N}=\{f:A\to\mathbb{R}^{n}:sup\{d(x,f(x)):x\in
A\}<1/30\}$ and $\mathcal{M}=\{f:\bar{B}^{n}\to\mathbb{R}^{n}\}$
then there exists a continuous map
$\phi:\mathcal{N}\to\mathcal{M}$ a constant $a_{0}$ and constants $b_{n}$ which depend only
on $n$ such that:
\begin{enumerate}
\item{$\phi(f)|_{S^{n-1}}=f|_{S^{n-1}}$ for all $f\in\mathcal{N}$}
\item{$\phi(id)=id$}
\item{If $f$ is $K$-quasiconformal, then $\phi(f)$ is
$a_{n}K^{3}$-quasiconformal}
\item{If f is locally L-bilipschitz, $\phi(f)$ is locally $a_{0}L^{3}$-bilipschitz}
\item{If f is L-bilipschitz, $\phi(f)$ is $a_{0}L^{3}$-bilipschitz}
\end{enumerate}
\end{thm}

\begin{rmk}
In the definition of $\mathcal{M}$ and $\mathcal{N}$ above we only
consider $f:A\to\mathbb{R}^{n}$ which are embeddings.
\end{rmk}

\newpage

\begin{lem}\label{lem:stab}
There exists $\epsilon>0$ so that
$$
\{f\in G(\mathbb{S}^{n}): \tilde{d}(id,f)<\epsilon\}\subset G(G_{0}(\mathbb{S}^{n}))
$$
\end{lem}

(Recall that $\tilde{d}(f,g)=sup\{d_{\mathbb{S}^{n}}(f(x),g(x)):x\in \mathbb{S}^{n}\}$).

\vspace{.3cm}

\noindent\textsc{Proof.} Take $f\in \mathrm{G}(\mathbb{S}^{n})$. Under the assumption that $f$ is close to $id$ we shall construct $g\in
\mathrm{G}_{0}(\mathbb{S}^{n})$ which agrees with $f$ on some open ball
$B\in \mathcal{B}$. This will mean that $g^{-1}\circ f$ will be
supported in the interior of $B^{c}$, and since
$int(B^{c})\in\mathcal{B}$, $g^{-1}\circ f$ will be in
$\mathrm{G}_{0}(\mathbb{S}^{n})$. As a result $f=g\circ (g^{-1}\circ f)$ will be an element of
$G(\mathrm{G}_{0}(\mathbb{S}^{n}))$ as required and the proof will be
complete.

\vspace{.3cm}

Begin by using Theorem \ref{thm:sul} to find a Sullivan group
$\Gamma$, giving rise to a Sullivan manifold
$Q=\mathbb{B}^{n}/\Gamma$ of dimension $n$. Let $\pi:\mathbb{B}^{n}\to Q$ denote the standard projection, then we have the following: For every $q\in Q$,
there exists $r(q)>0$ such that for any $p\in\pi^{-1}(q)$ there
exists a neighborhood $U$ of $\bar{B}_{\mathbb{B}^{n}}(r(q),p)$
such that the restriction $\pi|_{U}$ is an isometry. Since $Q$ is
compact we can fix $r>0$ which works for every $q\in Q$ and define
$D=B_{\mathbb{B}^{n}}(0,r)\subset \mathbb{B}^{n}$.

\vspace{.3cm}

Consider the restriction of $f$ to
the closed unit ball $\bar{\mathbb{B}}^{n}\subset \mathbb{S}^{n}$, which we
denote by $f_{0}$. We shall construct $g:\mathbb{S}^{n}\to\mathbb{S}^{n}$ which is supported on $\mathbb{B}^{n}$ and agrees with $f$ on $D$. It is thus necessary that we take $f$ sufficiently
close to the identity so as to ensure that $f_{0}(D)\subset
\mathbb{B}^{n}$. Choose $q\in Q\setminus \pi(\bar{D})$
and use condition 4 of Definition \ref{def:sul} to find a LIP
immersion $\iota:Q\setminus q\to B_{\mathbb{R}^{n}}(0,0.9)$.

\vspace{.3cm}

We can arrange that $\iota\circ\pi=id$ in a neighbourhood of $\bar{D}$. To see this, let $U$ be a neighbourhood of $\bar{D}$ on which $\iota\circ \pi$ is a LIP-embedding. Let $B_{1},B_{2}\subset \mathbb{B}^{n}$ denote open balls centred at $0$ which satisfy $\bar{D}\subset B_{1}\subset\bar{B}_{1}\subset B_{2}\subset U$. Apply Theorem \ref{thm:lip} to the map $(\iota\circ \pi)^{-1}$ on the annulus $B_{2}\setminus B_{1}$ to find a LIP-embedding $h':\mathbb{S}^{n}\setminus B_{1}\to\mathbb{S}^{n}$ which agrees with $(\iota\circ \pi)^{-1}$ on $\partial B_{1}$. Now define $h:\mathbb{S}^{n}\to \mathbb{S}^{n}$ by,
$$
h(x)=\left\{ \begin{array}{ll}
(\iota\circ \pi)^{-1}(x) & \textrm{if $x\in B_{1}$}\\
h'(x) & \textrm{if $x\in \mathbb{S}^{n}\setminus B_{1}$}\\
\end{array} \right.
$$
then $h$ is a lipeomorphism which agrees with $(\iota\circ\pi)^{-1}$ in a neighbourhood of $\bar{D}$. By post-composing $\iota$ with $h$ we will get a LIP immersion of $Q$ into $\mathbb{S}^{n}$ whose image is contained in some element of $\mathcal{B}$ which contains $\bar{D}$. We may now use a diffeomorphism which is the identity on some neighbourhood of $\bar{D}$ to map $h(B_{\mathbb{R}^{n}}(0,0.9))$ back inside $B_{\mathbb{R}^{n}}(0,0.9)$. This will give us a new LIP immersion $\iota':Q\setminus q\to B_{\mathbb{R}^{n}}(0,0.9)$ with the additional property that $\iota'\circ \pi=id$ in a neighbourhood of $\bar{D}$. We now assume that $\iota\circ\pi=id$ in a neighbourhood of $\bar{D}$.

\vspace{.3cm}

Choose $t>0$
sufficiently small so that $3t<r$ and for $i=1,2,3$ the balls
$D_{i}=B_{Q}(q,it)\subset Q$, satisfy
$D_{i}\cap\pi(\bar{D})=\phi$.

\newpage

We will construct $g$ by lifting maps as shown.

$$
\begin{array}{rccc}
g:             & \mathbb{S}^{n}        & \to & \mathbb{S}^{n} \\
               & \cup         &     &\cup \\
\tilde{f}_{2}: & \mathbb{B}^{n}        & \to & \mathbb{B}^{n} \\
               & \downarrow\pi&     & \downarrow\pi \\
f_{2}:         & Q        & \to & Q \\
               & \cup         &     & \cup\\
f_{1}:         & Q\setminus D_{2}&\to&Q\setminus D_{1}\\
               &\downarrow\iota&&\downarrow\iota\\
f_{0}:         &\bar{\mathbb{B}}^{n}&\to&\mathbb{S}^{n}\\
\end{array}
$$

Firstly we construct $f_{1}:Q\setminus D_{2}\to Q\setminus D_{1}$, this is a natural lift of $f_{0}$ constructed using the LIP-immersion $\iota$. In order to define $f_{1}$ we require that $f$ is
sufficiently close to $id$ (i.e. that $f_{0}$ is close to the
inclusion $\eta:\bar{B}^{n}\hookrightarrow\mathbb{R}^{n}\cup\{\infty\}$). This requirement is to ensure that $f_{0}\circ\iota(Q\setminus D_{2})$ is contained within the image of $\iota(Q\setminus D_{1})$, a necessary condition for the map $f_{1}$ to exist.

\vspace{.3cm}

To define $f_{1}$, we first construct a finite
open cover $U_{i}$ for $i=1,\ldots, k$ of $Q\setminus D_{2}$, such that the restriction
$\iota|_{U_{i}}$ is a LIP embedding for every $i$.

\vspace{.3cm}

Now construct another finite cover $W_{i}$ for $i=1,\dots,k$, this
time consisting of closed and hence compact subsets of $Q\setminus
D_{1}$ such that $W_{i}\subset U_{i}$. By taking $f$ closer to
$id$ if necessary we will have
$f_{0}(\iota(W_{i}))\subset\iota(U_{i})$ so we can define $f_{1}$
on each $W_{i}$ by
$f_{1}|_{W_{i}}=(\iota|_{U_{i}})^{-1}f_{0}(\iota|_{W_{i}})$. This
definition makes sense since $\iota$ is a LIP-embedding on each
$U_{i}$ and hence has a well defined inverse on its image.

\vspace{.3cm}

Now by choosing $f$ to lie in an even smaller neighbourhood of the
identity if necessary, we can apply the quantitative Schoenflies
theorem \ref{thm:quan} to the restriction of $f_{1}$ to $\bar{D}_{3}\setminus
D_{2}$ to construct a map $f_{2}:Q\to Q$ which agrees with
$f_{1}$ on $Q\setminus D_{3}$. Now we can lift $f_{2}$ to
$\tilde{f}_{2}:\mathbb{B}^{n}\to \mathbb{B}^{n}$ in the usual way so that
$\pi\circ \tilde{f}=f\circ \pi$. This lifting is unique if we specify
$\tilde{f}_{2}(0)\in\pi^{-1}\circ f_{2}\circ \pi(0)$, but we have only one
choice since we want $\tilde{f}_{2}|_{D}=f_{0}|_{D}$,
namely we set $f_{0}(0)=\tilde{f}_{2}(0)$.

\vspace{.3cm}

Now, if we choose $f$ to be still closer to the identity, then
since $\iota(Q\setminus q)$ is compactly contained in
$\mathbb{B}^{n}$, $f_{1}$ will move closer to the identity map of
$Q\setminus D_{2}$ (Here we are using the hyperbolic metric on $Q$
to determine this distance). Consequently the restriction of
$f_{1}$ to $\partial D_{3}$ will also become closer to the
identity (we used this to apply the quantitative Schoenflies
theorem \ref{thm:quan}). Since the map $\phi$ in theorem
\ref{thm:quan} is continuous and it fixes $id$,
$\phi(f_{1}|_{\bar{D}_{3}\setminus D_{2}})$ will also move closer
to the identity. As a result we are making $f_{2}$ close to the
identity in the hyperbolic metric on $Q$.

\vspace{.3cm}

If the dimension of $Q$ is 3 or more then by taking $f_{2}$ close
enough to the identity we can deduce that
$f_{2}$ is homotopic to the identity by applying Mostow rigidity.

\vspace{.3cm}

In dimension 2 however, we will show that $f_{2}$ will be homotopic to the identity by showing that it acts as the identity on $\pi_{1}$. To see this, take a set of generators $\tilde{\gamma}_{1}\cdots \tilde{\gamma}_{m}$ for $\pi_{1}$ (this will be a finite set as $Q$ is closed) and take a representative $\gamma_{i}$ of each. Then for each $\gamma_{i}$ there exists $\epsilon_{i}$ such that the $\epsilon_{i}$-neigbourhood of $\gamma_{i}$ is homotopy equivalent to $\gamma_{i}$. The same will then be true of $\epsilon'=$min$\{\epsilon_{i}:i=1,\dots,m\}$. One can now see that if $f_{2}$ is uniformly $\epsilon'$ close to the identity then $f_{2}(\gamma_{i})$ will be homotopic to $\gamma_{i}$ for each $i$. Consequently, it will induce the identity on $\pi_{1}$.

\vspace{.3cm}

The homotopy from $f_{2}$ to the identity then lifts to
show that $\tilde{f}_{2}$ also homotopic to the identity. Now,
since $Q$ is compact there exists a compact region
$K\subset\mathbb{B}^{n}$ such that $\pi(K)=Q$. Furthermore, the
lifted homotopy commutes with the action of $\Gamma$ and so is
completely determined by its behaviour on $K$. Now, each point in
$K$ can only be moved a finite hyperbolic distance by the
homotopy and hence the same applies to $\tilde{f}_{2}$. Since
$\Gamma$ acts by isometries this applies to every point in
$\mathbb{B}^{n}$ allowing us to deduce that $\tilde{f}_{2}$ is a
bounded hyperbolic distance from  the identity.

\vspace{.3cm}

We now use the fact that a lipeomorphism or K-quasiconformal
homeomorphism of hyperbolic space which is a bounded hyperbolic
distance from the identity will determine a map of the same class
from $\bar{B}^{n}$ to itself which is the identity on $\partial
\mathbb{B}^{n}$. This means that we can extend $\tilde{f}_{2}$ to
a map $g:\mathbb{S}^{n}\to \mathbb{S}^{n}$ which agrees with
$\tilde{f}_{2}$ on $\mathbb{B}^{n}$ and is the identity otherwise.

\vspace{.3cm}

From the construction of $g$ we have that $f_{0}\iota\pi=\iota\pi
g$ where defined (in particular this applies in a neighbourhood of
$D$) and $\iota\pi=id$ in a neighbourhood of $D$, consequently
$f=g$ on $D$ as we required.

\vspace{.3cm}

It remains to check that $g$ is in fact an element of
$G_{0}(\mathbb{S}^{n})$. Suppose that $f$ was quasiconformal and
let $E_{1}=\pi^{-1}(\bar{D}_{3})$ and $E_{2}=\pi^{-1}(Q\setminus
D_{3})$. Then $E_{1}$ consists of countably many components
$E_{1}^{i}$. Consider one such component $E_{1}^{1}$, then the dilatations of $\tilde{f}_{2}$ on $E_{1}^{1}$ are determined by the behaviour of $f_{2}$ on $\bar{D}_{3}$. Since $f_{2}|_{\bar{D}_{3}}$ was defined by
the quantitative Schoenflies theorem it is quasiconformal, and hence $\tilde{f}_{2}|E_{1}^{1}$ will be too.
Furthermore if $i\neq j$ there exists $\gamma\in\Gamma$ such that
$\tilde{f}_{2}|_{E_{1}^{i}}=\gamma\tilde{f}_{2}\gamma^{-1}|_{E_{1}^{j}}$
since $\gamma$ is conformal $g|_{E_{1}}$ is quasiconformal.

\vspace{.3cm}

Now since $Q$ is compact, there exists an open set $U\subset
\mathbb{B}^{n}$ whose closure in $\mathbb{B}^{n}$ is compact and such that
$\pi(U)=Q$, then $\iota\pi$ is locally bilipschitz on
$U\setminus\pi^{-1}(\bar{D}_{1})$. By construction,
$\iota\pi\tilde{f}_{2}=f_{0}\pi\iota$ in $E_{2}$ and hence
$\tilde{f}_{2}|_{E_{2}\cap U}$ is quasiconformal. Since
$\gamma\tilde{f}_{2}=\tilde{f}_{2}\gamma$ for all
$\gamma\in\Gamma$, we get $\tilde{f}_{2}|_{E_{2}}$ is quasiconformal. It now follows (see \cite{VA} Theorem 35.1) that $g$ is quasiconformal.

\vspace{.3cm}

If $f$ was a lipeomorphism of $\mathbb{S}^{n}$ then the same argument as above will allow us to deduce that $\tilde{f}_{2}$ is a lipeomorphism of $\mathbb{B}^{n}$ with its hyperbolic metric. This, along with the fact that $\tilde{f}_{2}$ is close to the identity, implies that $\tilde{f}_{2}$ is a lipeomorphism of $\mathbb{B}^{n}$ with the Euclidean metric and hence $g\in$ LIP($\mathbb{S}^{n}$).

\marginpar{$\square$}

\begin{thm}\label{thm:end}
$\mathrm{G}(\mathbb{S}^{n})=G(\mathrm{G}_{0}(\mathbb{S}^{n}))$
\end{thm}

\noindent\textsc{Proof.} $G(\mathrm{G}_{0}(\mathbb{S}^{n}))$ is a subgroup
of the topological group $\mathrm{G}(\mathbb{S}^{n})$ which by Lemma
\ref{lem:stab} contains a neighbourhood of the identity and
consequently we can translate this neighbourhood to show that
$G(\mathrm{G}_{0}(\mathbb{S}^{n}))$ is open in $\mathrm{G}(\mathbb{S}^{n})$. As
$G(\mathrm{G}_{0}(\mathbb{S}^{n}))$ is an open subgroup, it will be
closed, as its cosets will also be open. Consequently, it suffices
to show that $\mathrm{G}(\mathbb{S}^{n})$ is connected. So take $f\in
\mathrm{G}(\mathbb{S}^{n})$, we can easily construct a path in
$\mathrm{G}(\mathbb{S}^{n})$ which joins $f$ to an element with a fixed
point, so assume $f$ has a fixed point $x\in \mathbb{S}^{n}$. Stereographic
projection from $x$ allows us to identify $\mathbb{S}^{n}\setminus \{x\}$
with $\mathbb{R}^{n}$. We can consequently think of $f$ as a
G-homeomorphism of $\mathbb{R}^{n}$, in particular there will
exist $x_{0}\in\mathbb{R}^{n}$ at which $f$ is differentiable with
positive Jacobian. Furthermore, we may assume that
$x_{0}=0=f(x_{0})$. Set, $g_{t}(x)=f(tx)/t$ for $t>0$ and
$g_{0}=f'(0)$, then the path $g_{t}$ joins $f'(0)$ to $f$ and
since $det(f'(0))>0$ we can join $f'(0)$ to the identity. We can
extend all of these maps to $\mathbb{S}^{n}$ by leaving $x$ invariant
giving us a path in $\mathrm{G}(\mathbb{S}^{n})$ from $f$ to the identity.
This shows that $\mathrm{G}(\mathbb{S}^{n})$ is path connected and hence
connected. \marginpar{$\square$}

\begin{footnotesize}

\end{footnotesize}

\end{document}